\documentclass[11pt]{article}

\usepackage{latexsym,amsmath,amssymb,mathrsfs,graphicx,hyperref,mathpazo,xcolor,amsthm}
\hypersetup{colorlinks=true,linkcolor=blue,citecolor=blue,filecolor=blue,urlcolor=blue,
  pdftitle={A systematic approach to Diophantine equations: open problems},
  pdfauthor={Bogdan Grechuk}}

\usepackage{url}
\definecolor{reviewgreen}{rgb}{0,0.45,0}

\def\tr{{\raise0pt\hbox{$\scriptscriptstyle\top$}}}

\newtheorem{theorem}{Theorem}

\newtheorem{problem}[theorem]{Problem}
\newtheorem{definition}{Definition}

\title{\vspace{-1.5 cm}{\bf A systematic approach to Diophantine equations: open problems}}

\author{Bogdan Grechuk}

\begin{document}

\maketitle

\begin{abstract}
	This paper collects polynomial Diophantine equations that are simple to state but apparently difficult to solve. 
\end{abstract}

\section{Background and definitions}

This document collects the smallest Diophantine equations that are currently open, in the sense described in \cite{mainbook}. We plan to update this document regularly to keep the list current.

We start by recalling the necessary definitions from \cite[Section 1.1.1]{mainbook}. A monomial in variables $x_1,\dots,x_n$ with an integer coefficient is any expression of the form $M=ax_1^{k_1}\dots x_n^{k_n}$, where $a$ is a nonzero integer and $k_1, \dots, k_n$ are nonnegative integers. The integer $d=\sum_{i=1}^n k_i$ is called the degree of the monomial $M$. Two monomials $M=ax_1^{k_1}\dots x_n^{k_n}$ and $M'=a'x_1^{k'_1}\dots x_n^{k'_n}$ are called similar if $k_i=k'_i$ for all $i=1,\dots,n$. 

A polynomial $P(x_1,\dots,x_n)$ with integer coefficients is a finite sum of monomials. A polynomial Diophantine equation is an equation of the form 
\begin{equation}\label{eq:diofgen}
	P(x_1, \dots, x_n) = 0.
\end{equation} 
We assume that at least one monomial of $P$ is non-constant, so that equation \eqref{eq:diofgen} has at least one variable. Variables that do not occur in $P$ are omitted.
Without loss of generality, we may also assume that $P$ is in reduced form; that is, no two monomials of $P$ are similar. If the polynomial $P$ consists of $k$ non-similar monomials with integer coefficients $a_1, \dots, a_k$ and degrees $d_1, \dots, d_k$, respectively, then the size of equation \eqref{eq:diofgen} is defined as
\begin{equation}\label{eq:Hdef}
	H(P)=\sum_{i=1}^k |a_i| 2^{d_i}.
\end{equation}
To compute the size $H(P)$ of any equation $P = 0$, we simply substitute $2$ for each variable, replace each coefficient by its absolute value, and evaluate. It follows directly from the definition that, up to renaming the variables, there are only finitely many polynomial Diophantine equations of any given size. Thus, we can order the equations by their size $H$ and study them in that order.

The aim of this document is to list the smallest equations that are currently open. We will list only one equation from each equivalence class. For unrestricted integer solutions, two equations of the form \eqref{eq:diofgen} are called equivalent if one can be transformed into the other by a sequence of the following operations: (i) multiplying the defining polynomial by $-1$; (ii) replacing any variable $x_i$ by $-x_i$; and (iii) renaming or permuting the variables. For positive integer solutions, only operations (i) and (iii) are used, since operation (ii) does not preserve positivity.

We will consider unrestricted equations as well as equations in various categories. For example, we may restrict the number of variables in equation \eqref{eq:diofgen}, its degree (defined as the maximum degree of a monomial of $P$), or the number of monomials in $P$. We will also consider the following categories.
\begin{itemize}
	\item \textbf{Homogeneous equations}: equations in which all monomials have the same degree.
	\item \textbf{Symmetric equations}: equations invariant under permutations of the variables. In other words, if a symmetric equation contains a monomial $ax_1^{\alpha_1}\dots x_n^{\alpha_n}$ and $(\beta_1, \dots, \beta_n)$ is any permutation of $(\alpha_1,\dots,\alpha_n)$, then the equation must also contain a monomial $ax_1^{\beta_1}\dots x_n^{\beta_n}$.
	\item \textbf{Cyclic equations}: equations invariant under the cyclic shift of variables $(x_1,\dots,x_n)\mapsto(x_2,\dots,x_n,x_1)$. Equivalently, if the defining polynomial contains $a\prod_{i=1}^n x_i^{\alpha_i}$, it contains
	\[
	a\prod_{i=1}^n x_{i+j}^{\alpha_i},\qquad j=0,\dots,n-1,
	\]
	where variable indices are interpreted modulo $n$.
	\item \textbf{Equations with independent monomials}: equations in which no two monomials share a variable. 
\end{itemize}
 
One may ask various questions about equation \eqref{eq:diofgen}. For example, one may ask how to describe all its integer (or rational) solutions, whether its solution set is finite, or whether it is nonempty. These questions will be discussed in the following sections.

The claims that the equations listed below are the smallest or shortest open equations in their respective categories were also cross-checked with the assistance of ChatGPT. Under the author's guidance, ChatGPT produced explicit solutions to the corresponding problems for every other equation up to the stated size or length bound.

Version 1 of this document corresponds exactly to the equations left open in \cite{mainbook}. In subsequent versions, we have removed equations that were solved and added equations that became the smallest/shortest open examples after earlier problems were resolved. All new equations that are not mentioned in \cite{mainbook} are \textcolor{red}{marked in red}.

\section{Polynomial parametrization of all integer solutions}

One natural question about \eqref{eq:diofgen} is whether all its integer solutions can be parametrized.

\begin{definition}\label{def:polfam}
	Let ${\mathbb Z}^n$ be the set of vectors $(x_1,\dots, x_n)$ with integer coordinates $x_i$. We say that a subset $S \subseteq {\mathbb Z}^n$ is a polynomial family if there exist an integer $k\geq 0$ and polynomials $P_1, \dots, P_n$ in $k$ variables $u_1,\dots, u_k$ with integer coefficients such that $(x_1,\dots, x_n)\in S$ if and only if there exist integers $u_1,\dots, u_k$ such that 
	$$
	x_i=P_i(u_1,\dots, u_k), \quad i=1,\dots, n.
	$$ 
\end{definition}

\begin{problem}\label{prob:H12existpol}
	Given a polynomial Diophantine equation \eqref{eq:diofgen}, determine whether its set of integer solutions is a finite union of polynomial families. 
\end{problem} 

The smallest equations for which Problem \ref{prob:H12existpol} remains open are listed in Table \ref{tab:H13open}.

\begin{table}
	\begin{center}
		\begin{tabular}{ |c|c|c|c| } 
			\hline
			$H$ & Equation & $H$ & Equation \\ 
			\hline\hline
			$13$ & $xyz+t^2+1=0$  & $13$ & $x^3+yz+1=0$ \\ 
			\hline
		\end{tabular}
		\caption{\label{tab:H13open} Equations of size $H \leq 13$ for which the existence of a polynomial parametrization is unknown.}
	\end{center} 
\end{table}

\section{Describing all integer solutions}\label{sec:allsol}

In addition to polynomial parametrization, there are many other ways to describe all integer solutions of \eqref{eq:diofgen}. For example, one may use
\begin{itemize}
	\item rational expressions involving a parameter in the denominator, together with the restriction that the parameter divide the numerator, thereby ensuring that the expressions produce integer solutions;
	\item recurrence relations; or
	\item a set of ``initial'' solutions and a set of ``transformations'' that map solutions to solutions, such that all solutions of \eqref{eq:diofgen} can be obtained from the initial solutions by applying a sequence of transformations.
\end{itemize} 

Given a polynomial Diophantine equation \eqref{eq:diofgen}, one may ask whether its set of integer solutions can be described in any ``reasonable'' way. Because ``reasonable'' is not defined rigorously here, this question is necessarily informal. At present, we regard the solution sets of all equations of size $H\leq 17$ as having reasonably explicit descriptions. In the range $18\leq H\leq 21$, some equations have only complicated known descriptions, so their classification as ``open'' or ``solved'' is debatable. For example, to describe all integer solutions to equations 
\begin{equation}\label{eq:zy2pm1x2p2}
	\textcolor{red}{z(y^2\pm 1)=x^3+2}
\end{equation}
of size $H=20$, define, for each nonzero integer $n$,
\[
S_{\pm}(n)=\{v\in\mathbb Z: v^2\pm1\neq0\ \text{and}\ v^2\pm1\mid n\}.
\]
Each set $S_{\pm}(n)$ is finite, contains $0$, and can be computed by enumerating the signed divisors of $n$. Since $u^3+2\neq0$ for every integer $u$, all integer solutions to \eqref{eq:zy2pm1x2p2} are exactly
\[
(x,y,z)=\left(u,v,\frac{u^3+2}{v^2\pm1}\right),
\qquad u\in\mathbb Z,\quad v\in S_{\pm}(u^3+2).
\]
The excluded values $v=\pm1$ in the minus case cannot give solutions, since they would require $u^3+2=0$. It is debatable whether this is an acceptable description of all integer solutions to \eqref{eq:zy2pm1x2p2}; we currently have no description in more conventional forms, such as polynomial parametrizations or recurrence relations.

The smallest equation that is unambiguously open is
\begin{equation}\label{eq:z2py2zpx3m2}
	z^2+y^2z+x^3-2=0
\end{equation}	
of size $H=22$. For this equation, it is not even known whether it has finitely or infinitely many integer solutions.

We next consider equations in various categories. 
The smallest open cyclic equation is 
\begin{equation}\label{eq:x2ypy2zpz2xm1}
	x^2y + y^2z + z^2x = 1
\end{equation}
of size $H=25$. For this equation, we do not know whether its integer solution set is finite or infinite. 
The smallest open symmetric equation is 
\begin{equation}\label{eq:x3py3pz3m2} 
	x^3+y^3+z^3=2
\end{equation}
of size $H=26$. For this equation, we do not even know whether there are infinitely many integer solutions other than permutations of those in the family $(x,y,z)=(1-6t^3, -6t^2, 1+6t^3)$, $t \in {\mathbb Z}$. 

\begin{table}
	\begin{center}
		\begin{tabular}{ |c|c|c|c| } 
			\hline
			$H$ & Equation & $H$ & Equation \\ 
			\hline\hline
			$29$ & $x^4+y^3+xy-1=0$ & $31$ & $x^4-y^3+xy+x+1=0$ \\ 
			\hline
			$30$ & $x^4+y^3+x+2y=0$ & $31$ & $x^4+y^3+xy-3=0$ \\ 
			\hline
			$30$ & $x^4+y^3+x-y-2=0$ & $31$ & $x^4+y^3+x-2y-1=0$ \\ 
			\hline
			$30$ & $x^4+y^3+x+y+2=0$ & $31$ & $x^4+y^3+xy+x-1=0$ \\ 
			\hline
			$30$ & $x^4+y^3+x+y-2=0$ & $31$ & $x^4+y^3+x+2y+1=0$ \\ 
			\hline
		\end{tabular}
		\caption{\label{tab:H31Open} Open two-variable equations of size $H\leq 31$.}
	\end{center} 
\end{table}

The smallest open three-monomial equations are
\begin{equation}\label{eq:yx3myz2mz}
	(a) \,\, y(x^3-z^2)=z \quad \quad (b) \,\, x^2y^2+x=z^3
\end{equation}
and
\begin{equation}\label{eq:yx3myz2mx}
	y(x^3-z^2)=x
\end{equation}
of size $H=26$. Equations \eqref{eq:yx3myz2mz}(a) and \eqref{eq:yx3myz2mz}(b) are mutually reducible, and each reduces further to
\begin{equation}\label{eq:x3y2ez3m1}
	x^3y^2=z^3-1,
\end{equation}
as explained in \cite[Section 4.3.5]{mainbook}. Similarly, equation \eqref{eq:yx3myz2mx} reduces to the equations $x^4y^3=z^2\pm 1$. These equations remain open: some obvious solutions are known (see \cite[Section 4.3.5]{mainbook}), but it is not known whether any others exist.

The smallest open two-variable equations are listed in Table \ref{tab:H31Open}. For each of these equations, it is known that the integer solution set is finite, and the remaining open problem is to list all integer solutions. 

We may also consider various intersections of the above categories. For example, the smallest open symmetric homogeneous equation is
\begin{equation}\label{eq:x3py3pz3pt3ps3}
	x^3+y^3+z^3+t^3+s^3=0
\end{equation} 
of size $H=40$. The smallest open cyclic homogeneous equations are \eqref{eq:x3py3pz3pt3ps3} and the equation
\begin{equation}\label{eq:x2ypy2zpz2tpt2sps2x}
	x^2y+y^2z+z^2t+t^2s+s^2x=0
\end{equation} 
of the same size. The smallest open symmetric two-variable equations are listed in Table \ref{H325sym2varopen}.

\begin{table}
	\begin{center}
		\begin{tabular}{ |c|c| } 
			\hline
			$H$ & Equation \\ 
			\hline\hline
			$325$ & $x^6 + x^5 y + x^3 y^3 + x y^5 + y^6 + x + y + 1=0$ \\ 
			\hline
			$325$ & $x^6 - x^5 y - x^3 y^3 - x y^5 + y^6 + x + y - 1=0$ \\ 
			\hline
			$325$ & $x^6 - x^4 y^2 + x^3 y^3 - x^2 y^4 + y^6 + x + y + 1=0$ \\ 
			\hline
			$325$ & $x^6 - x^4 y^2 - x^3 y^3 - x^2 y^4 + y^6 + x + y - 1=0$ \\ 
			\hline
			$325$ & $x^6 + 3 x^3 y^3 + y^6 + x + y + 1=0$ \\ 
			\hline
			$325$ & $x^6 - 3 x^3 y^3 + y^6 + x + y - 1=0$ \\ 
			\hline
		\end{tabular}
		\caption{\label{H325sym2varopen} The smallest open symmetric two-variable equations.}
	\end{center} 
\end{table}

\section{Existence of arbitrarily large integer solutions}\label{sec:large}

The question of describing all integer solutions in a ``reasonable'' way is informal. This motivates us to consider the following restricted but more formal question.

\begin{problem}\label{prob:large}
	Given a Diophantine equation \eqref{eq:diofgen}, determine whether, for every $k\geq 0$, it has a solution satisfying
	\begin{equation}\label{eq:largecond}
		\min(|x_1|, \dots, |x_n|) \geq k.
	\end{equation}
	If so, the problem is solved. If not, describe all its integer solutions.
\end{problem}

The smallest equations for which Problem \ref{prob:large} is open are listed in Table \ref{tab:H23openlarge}.

\begin{table}
	\begin{center}
		\begin{tabular}{ |c|c|c|c| } 
			\hline
			$H$ & Equation &  $H$ & Equation \\ 
			\hline\hline
			$22$ & $z^2+y^2z+x^3-2=0$ &  $23$ & $z^2+y^2z+x^3-x-1=0$ \\ 
			\hline
			$23$ & $z^2+y^2z+x^3-3=0$ &  $23$ & \textcolor{red}{$y^2+x^2y+z^2x+y-1=0$} \\ 
			\hline
			$23$ & $z^2+y^2z+x^3+3=0$ &   &  \\ 
			\hline
		\end{tabular}
		\caption{\label{tab:H23openlarge} Equations of size $H\leq 23$ for which Problem \ref{prob:large} is open.}
	\end{center} 
\end{table}

We can also consider equations in various categories. The smallest cyclic equation for which Problem \ref{prob:large} remains open is \eqref{eq:x2ypy2zpz2xm1}. The smallest open symmetric equation is
\begin{equation}\label{eq:x3py3pz3m3}
	x^3+y^3+z^3 = 3
\end{equation}
of size $H=27$. This is also the smallest open equation with independent monomials. The next-smallest such equation is
\begin{equation}\label{eq:x4py3pz2p1}
	x^4+y^3+z^2+1=0
\end{equation} 
of size $H=29$. The smallest open two-variable equations are those listed in Table \ref{tab:H31Open}. The smallest open three-monomial equations are \eqref{eq:yx3myz2mz} and \eqref{eq:yx3myz2mx}.

For a homogeneous equation, the existence of solutions satisfying \eqref{eq:largecond} for every $k$ is equivalent, by scaling, to the following condition. A negative answer still leaves the problem of describing solutions with a zero coordinate.

\begin{problem}\label{prob:homlarge}
	Given a homogeneous polynomial $P(x_1,\dots,x_n)$ with integer coefficients, determine whether equation \eqref{eq:diofgen} has an integer solution in which all variables are nonzero.
\end{problem}

The smallest homogeneous equations for which Problem \ref{prob:homlarge} is open are listed in Table \ref{tab:H64genus3}.

\begin{table}
	\begin{center}
		\begin{tabular}{ |c|c|c|c| } 
			\hline
			$H$ & Equation & $H$ & Equation \\ 
			\hline\hline
			$64$ & $x^4+x y^3+z^4+t^4 = 0$ & $64$ & $x^4+y^4+z^3t-zt^3 = 0$ \\ 
			\hline			
		\end{tabular}
		\caption{\label{tab:H64genus3} Equations of size $H\leq 64$ for which Problem \ref{prob:homlarge} is open.}
	\end{center} 
\end{table}

One may also study Problem \ref{prob:homlarge} for diagonal homogeneous equations of the form 
\begin{equation}\label{eq:homdiand}
	a_1x_1^d+a_2x_2^d+\dots + a_nx_n^d=0.
\end{equation}
The smallest equation of the form \eqref{eq:homdiand} for which Problem \ref{prob:homlarge} remains open is
\begin{equation}\label{eq:3x4py4mz4mt4}
	\textcolor{red}{3x^4+y^4=z^4+t^4}
\end{equation} 
of size $H=96$.

\section{The finiteness problem}\label{sec:fin}

In general, Problem \ref{prob:large} is not fully formal because the word ``describe'' in its second part is ambiguous. We therefore consider the following narrower, fully rigorous question.

\begin{problem}\label{prob:fin}
	Given a Diophantine equation, either list all its integer solutions or prove that there are infinitely many of them. 
\end{problem}

The smallest equations for which Problem \ref{prob:fin} is open are listed in Table \ref{tab:H24finopen}.

\begin{table}
	\begin{center}
		\begin{tabular}{ |c|c|c|c| } 
			\hline
			$H$ & Equation &  $H$ & Equation \\ 
			\hline\hline
			$22$ & $z^2+y^2z+x^3-2=0$ &  $23$ & $z^2+y^2z+x^3-x-1=0$ \\ 
			\hline
			$23$ & $z^2+y^2z+x^3-3=0$ &  $24$ & $z^2+y^2z+x^3-x-2=0$  \\ 
			\hline
			$23$ & $z^2+y^2z+x^3+3=0$ &  $24$ & $z^2+y^2z+x^3-x+2=0$  \\ 
			\hline
		\end{tabular}
		\caption{\label{tab:H24finopen} Equations of size $H\leq 24$ for which Problem \ref{prob:fin} is open.}
	\end{center} 
\end{table}

The smallest open cyclic and symmetric equations are \eqref{eq:x2ypy2zpz2xm1} and \eqref{eq:x3py3pz3m3}, respectively. The smallest and next-smallest equations with independent monomials are \eqref{eq:x3py3pz3m3} and \eqref{eq:x4py3pz2p1}. The smallest open two-variable equations are those listed in Table \ref{tab:H31Open}. The smallest open three-monomial equation is
\begin{equation}\label{eq:x3y2mz3m2}
	x^3y^2 = z^3+2
\end{equation}
of size $H=42$. For homogeneous equations, Problem \ref{prob:fin} reduces to the following question.

\begin{problem}\label{prob:finhom}
	Given a homogeneous polynomial $P$ with integer coefficients, determine whether equation \eqref{eq:diofgen} has an integer solution $(x_1, \dots, x_n) \neq (0, \dots, 0)$.
\end{problem}

The smallest homogeneous equation for which Problem \ref{prob:finhom} is open is
\begin{equation}\label{eq:x4px3ymy4py3zpz4}
	x^4+x^3 y-y^4+y^3 z+z^4=0
\end{equation}
of size $H=80$. This is the only such open equation of size $H\leq 100$.

The smallest equations of the form \eqref{eq:homdiand} for which Problem \ref{prob:finhom} remains open are listed in Table \ref{tab:H288hom44}.

\begin{table}
	\begin{center}
		\begin{tabular}{ |c|c|c|c| } 
			\hline
			$H$ & Equation &  $H$ & Equation \\ 
			\hline\hline
			$288$ & $7x^4+6y^4+4z^4-t^4= 0$ & $288$ & $7x^4+6y^4+3z^4-2t^4= 0$ \\ 
			\hline	
		\end{tabular}
		\caption{\label{tab:H288hom44} Open equations of the form \eqref{eq:homdiand} and size $H\leq 288$.}
	\end{center} 
\end{table}

\section{Existence of integer solutions}\label{sec:existence}

Perhaps the most famous question about Diophantine equations is Hilbert's tenth problem:

\begin{problem}\label{prob:main}
	Given a Diophantine equation, determine whether it has an integer solution.
\end{problem}

The smallest equations for which Problem \ref{prob:main} is open are listed in Table \ref{tab:H34unsolved}. 

\begin{table}
	\begin{center}
		\begin{tabular}{ |c|c|c|c| } 
			\hline
			$H$ & Equation & $H$ & Equation \\ 
			\hline\hline
			 $34$ & $6+x^3-y^2+y^2 z^2=0$ & $34$ &  $2+x^3+x y+y^2+y^2 z^2=0$ \\ 
			\hline
			$34$ & $6+x^3+x y^2 z+z^2=0$ &  &  \\ 
			\hline	
		\end{tabular}
		\caption{\label{tab:H34unsolved} Equations of size $H\leq 34$ for which Problem \ref{prob:main} is open.}
	\end{center} 
\end{table}

The smallest open equations with independent monomials are
\begin{equation}\label{eq:x4py3pz3p4}
	x^4 + y^3 + z^3 = -4
\end{equation}
and
\begin{equation}\label{eq:x4py3pz3m4}
	x^4 + y^3 + z^3 = 4
\end{equation} 
of size $H=36$. The smallest open cubic equation is 
\begin{equation}\label{eq:H39sym}
	x^3+x+y^3+y+z^3+z = xyz + 1
\end{equation}
of size $H=39$. This equation is also the smallest open symmetric equation. The next-smallest open symmetric equation is
$$
x^3-x+y^3-y+z^3-z=xyz+7
$$
of size $H=45$. The same two symmetric equations are also the smallest and next-smallest open cyclic equations. The smallest three-monomial equation for which Problem \ref{prob:main} is open is
\begin{equation}\label{eq:H46term3}
	x^3y^2 = z^3 + 6
\end{equation}
of size $H=46$.

One may also ask whether solutions exist in positive integers.

\begin{problem}\label{prob:pos}
	Given a Diophantine equation, determine whether it has a solution in \textbf{positive} integers.
\end{problem}

The smallest equations for which Problem \ref{prob:pos} remains open are the three equations \eqref{eq:yx3myz2mz}(a), \eqref{eq:yx3myz2mz}(b), and \eqref{eq:yx3myz2mx}; each has size $H=26$. As noted above, equations \eqref{eq:yx3myz2mz}(a) and \eqref{eq:yx3myz2mz}(b) are mutually reducible.

\section{The shortest open equations}\label{sec:shortopen}

There are other natural ways to order polynomial Diophantine equations. One is to define the \textbf{length} of a polynomial $P$ consisting of monomials of degrees $d_1, \dots, d_k$ with coefficients $a_1, \dots, a_k$ by
\begin{equation}\label{eq:ldef}
	l(P) = \sum_{i=1}^k  \log_2(|a_i|) + \sum_{i=1}^k d_i = \log_2 L(P), \quad \text{where} \quad L(P) = \prod_{i=1}^k\left(|a_i|\cdot 2^{d_i}\right).
\end{equation}
Ordering polynomials by length $l(P)$ is equivalent to ordering them by $L(P)$.

The equations in Table \ref{tab:H13open} have length $l=5$, and they are the shortest equations for which Problem \ref{prob:H12existpol} remains open. If descriptions other than polynomial parametrizations are allowed, every equation of length $l<8$ has a reasonably explicit description of its solution set, whereas the family of length-$8$ equations contains examples like \eqref{eq:zy2pm1x2p2} with uncertain status as well as equations in Table \ref{tab:l8openlarge} that are unambiguously open. The shortest open two-variable equation is $x^4+xy+y^3-1=0$, listed in Table \ref{tab:l9openfin}; the shortest open symmetric two-variable equations are the last two equations in Table \ref{H325sym2varopen}. 

\begin{table}[htbp]
	\begin{center}
		\begin{tabular}{ |c|c|c|c| } 
			\hline
			$l$ & Equation & $l$ & Equation \\ 
			\hline\hline
			$8$ & $y(x^3-z^2)=x$ & $8$ & $y(x^3-z^2)=z \,\, \Leftrightarrow \,\, x^2y^2+x=z^3 \,\, \Leftrightarrow \,\, x^3y^2=z^3-1$  \\ 
			\hline
			$8$ & $y(x^3-z^2)=x+1$ &  &  \\ 
			\hline
		\end{tabular}
		\caption{\label{tab:l8openlarge} Equations of length $l\leq 8$ for which Problem \ref{prob:large} is open.}
	\end{center} 
\end{table} 

\begin{table}[htbp]
	\begin{center}
		\begin{tabular}{ |c|c|c|c| } 
			\hline
			$l$ & Equation &  $l$ & Equation \\ 
			\hline\hline
			$9$ & $z^2+y^2z+x^3-2=0$ &  $9$ & $x^2y+y^2z+z^2x=1$ \\ 
			\hline
			$9$ & $z^2+y^2z+x^3-x-1=0$ &  $9$ & \textcolor{red}{$y^2+x^2y+z^2x+y-1=0$} \\ 
			\hline
		\end{tabular}
		\caption{\label{tab:l9openlargecub} Cubic equations of length $l\leq 9$ for which Problem \ref{prob:large} is open.}
	\end{center} 
\end{table}

\begin{table}
	\begin{center}
		\begin{tabular}{ |c|c|c|c|c|c| } 
			\hline
			$l$ & Equation & $l$ & Equation & $l$ & Equation \\ 
			\hline\hline
			$9$ & $z^2+y^2z+x^3-2=0$ & $9$ & $y(x^3-z^2)=x^2+1$ & $9$ & $x^3y^2=z^3+z+1$ \\ 
			\hline
			$9$ & $z^2+y^2z+x^3-x-1=0$ & $9$ & $y(x^3-z^2)=2x-1$ & $9$ & $x^3y^2=2z^3+1$ \\ 
			\hline
			$9$ & $x^2y+y^2z+z^2x=1$ & $9$ & $y(x^3-z^2)=2x+1$ & $9$ & $x^3y^2=z^4+1$ \\ 
			\hline
			$9$ & $z^2+y^2z+x^3y+1=0$ & $9$ & $x^3y^2=z^3+2$ & $9$ & $x^4y^3=z^2+1$ \\ 
			\hline
			$9$ & $x^4+y^3+z^2+1=0$ & $9$ & $x^3y^2=z^3-z+1$ & $9$ & $x^4+xy+y^3-1=0$ \\ 
			\hline
		\end{tabular}
		\caption{\label{tab:l9openfin} Equations of length $l\leq 9$ for which Problem \ref{prob:fin} is open.}
	\end{center} 
\end{table} 

The shortest equations for which Problems \ref{prob:large}, \ref{prob:fin}, and \ref{prob:main} are open are listed in Tables \ref{tab:l8openlarge}, \ref{tab:l9openfin}, and \ref{tab:lless11openexist}, respectively.

All equations in Tables \ref{tab:l8openlarge} and \ref{tab:lless11openexist} have degree at least $4$. The shortest cubic equations for which Problems \ref{prob:large} and \ref{prob:main} remain open are listed in Tables \ref{tab:l9openlargecub} and \ref{tab:cubshortest}, respectively. 

\begin{table}
	\begin{center}
		\begin{tabular}{ |c|c|c|c|c|c| } 
			\hline
			$l$ & Equation & $l$ & Equation & $l$ & Equation \\ 
			\hline\hline
			$10$ & $y^2+x^3y+z^4+1=0$ & $\approx 10.6$ & $x^3y^2=z^3+6$ & $\approx 10.6$ & $x^4y^3=z^2+3$  \\ 
			\hline
			$10$ & $x^3y^2=z^4+2$ & $\approx 10.6$ &  $x^3y^2=z^4-3$  &  & \\ 
			\hline
		\end{tabular}
		\caption{\label{tab:lless11openexist} Equations of length $l<11$ for which Problem \ref{prob:main} is open.}
	\end{center} 
\end{table}

For both Problems \ref{prob:large} and \ref{prob:fin}, the shortest open cyclic and symmetric equations are \eqref{eq:x2ypy2zpz2xm1} and \eqref{eq:x3py3pz3m3}, respectively. For Problem \ref{prob:main}, the shortest open symmetric equation is \eqref{eq:H39sym}, whose length is $l=15$. The next-shortest is the famous open equation
\begin{equation}\label{eq:x3py3pz3e114}
	\textcolor{red}{x^3+y^3+z^3 = 114},
\end{equation} 
whose length is $9+\log_2 114 \approx 15.83$. These two equations are also the shortest and next-shortest \emph{cyclic} equations for which Problem \ref{prob:main} remains open.

\begin{table}
	\begin{center}
		\begin{tabular}{ |c|c|c|c| }
			\hline
			$l$ & Equation & $l$ & Equation \\ 
			\hline\hline
			$\approx 13.6$ &  \textcolor{red}{$(3x-1)y^2 + x z^2 = x^3-2$} & $\approx 13.6$ &  \textcolor{red}{$y^2z+yz^2=6x^3+x^2+1$} \\
			\hline
		\end{tabular}
		\caption{\label{tab:cubshortest} The shortest cubic equations for which Problem \ref{prob:main} is open.}
	\end{center}
\end{table}

The shortest homogeneous equations for which Problem \ref{prob:homlarge} remains open are those in Table \ref{tab:H64genus3}. The shortest equation of the form \eqref{eq:homdiand} for which Problem \ref{prob:homlarge} remains open is \eqref{eq:3x4py4mz4mt4}. The shortest homogeneous equation for which Problem \ref{prob:finhom} remains open is \eqref{eq:x4px3ymy4py3zpz4}, whose length is $l=20$. The shortest equation of the form \eqref{eq:homdiand} for which Problem \ref{prob:finhom} remains open is
\begin{equation}\label{eq:x4py4pz4e51t4}
	\textcolor{red}{x^4+y^4+z^4=51t^4}
\end{equation} 
of length $l=16+\log_2 51 \approx 21.67$. Problem \ref{prob:finhom} for this equation is equivalent to asking whether $51$ is a sum of at most three rational fourth powers.\footnote{See \url{https://mathoverflow.net/questions/514531/}}

The shortest equations for which Problem \ref{prob:pos} remains open are \eqref{eq:yx3myz2mz}(a), \eqref{eq:yx3myz2mz}(b), \eqref{eq:yx3myz2mx}, and \eqref{eq:x3y2ez3m1}, all of length $l=8$. As explained above, equations \eqref{eq:yx3myz2mz}(a), \eqref{eq:yx3myz2mz}(b), and \eqref{eq:x3y2ez3m1} are mutually reducible. Thus, up to these reductions, only two open equations of length $l\leq 8$ remain: \eqref{eq:yx3myz2mx} and \eqref{eq:x3y2ez3m1}.

\section{Changes in this document between versions}

\subsection{Changes between versions 1 and 2}

Version 1 of this document corresponds exactly to the equations left open in \cite{mainbook}. Below is a summary of the updates between versions 1 and 2.

While preparing the survey \cite{ratcliffe2025generalized}, Ashleigh Ratcliffe noted that the primitive integer solutions to $x^4+y^4+2z^3=0$ had been determined in \cite[Theorem 8.1]{bennett2004}, and those to $x^5+y^4+2z^2=0$ in \cite[Theorem 1]{bennett2010}. These equations have therefore been excluded from Section \ref{sec:allsol}.

For the equations $2x^2-xyz-y^2-1=0$, $2x^2-xyz+y^2+1=0$, and $y^2+xyz+x^3+2=0$, Problem \ref{prob:large} has been solved\footnote{See \url{https://mathoverflow.net/questions/449559/} and \url{https://mathoverflow.net/questions/411958/}}. These equations have therefore been removed from Tables \ref{tab:H23openlarge} and \ref{tab:l8openlarge}.

For the equations $y^2+xyz+x^3+2=0$, $y^2+xyz+x^3-3=0$, $y^2+xyz+x^3+4=0$, $y^2+xyz+x^3-x-2=0$, $y^2+xyz+y+x^3+2=0$, and $y^2+xyz+2x^3+1=0$, Problem \ref{prob:fin} has been solved\footnote{See \url{https://mathoverflow.net/questions/411958/} and \url{https://mathoverflow.net/questions/466803/}}. These equations have therefore been removed from Tables \ref{tab:H24finopen} and \ref{tab:l9openfin}.

Problem \ref{prob:main} has been solved for the equations $6x^2z+z^2y+3y^3+1=0$, $y^2+10xyz+x^3-x-2=0$, and $y^2+7xyz+3x^3-2=0$.\footnote{See \url{https://mathoverflow.net/questions/466803/}} These were previously the shortest cubic equations for which Problem \ref{prob:main} was open, so they have been removed from Section \ref{sec:shortopen}.

\subsection{Changes between versions 2 and 3}

Ashleigh Ratcliffe reported that $2x^4+y^3+z^3=0$ is solved in \cite{zhang2014}; see also \cite[Theorem 4.3]{ratcliffe2025generalized}. She also determined the primitive integer solutions to $x^4+2y^3+2z^3=0$, $x^4-y^4+2z^3=0$, and $x^4-y^4+3z^3=0$; see \cite[Propositions A.2--A.4]{ratcliffe2024}. These equations have therefore been excluded from Section \ref{sec:allsol}. Furthermore, the equation $2x^4+2y^3+z^3=0$ has been removed from Section \ref{sec:allsol} because it duplicates $x^4+4y^3+z^3=0$: indeed, its solutions necessarily have even $z$, and substituting $z=2u$ and dividing by $2$ gives $x^4+y^3+4u^3=0$. Renaming $u$ as $y$ and the original $y$ as $z$ gives $x^4+4y^3+z^3=0$.

\subsection{Changes between versions 3 and 4}

Problem \ref{prob:fin} has been solved for $x^2y^2+x^2z-z^2-1=0$. The substitution $z\mapsto xy+(x^2-u)/2$, where $u$ is a new variable, transforms this equation into $u^2-4uxy-x^4+4=0$, which can be handled by the method used for the equations removed between versions 1 and 2, see \cite{grechuk2026cubic}. Jeremy Rouse also solved Problem \ref{prob:fin} for $y^2+z^2=x^5-1$ by a simple argument that applies equally to $y^2+2x^3y+z^2+1=0$, see \cite{grechuk2026cubic}. These three equations have therefore been removed from Table \ref{tab:l9openfin}.

Nguyen Xuan Tho proved that $7x^4-7y^4=25z^4$ has no integer solutions with $xyz\neq0$ \cite[Theorem 1.1]{tho2026} and has therefore been removed from Section \ref{sec:large}.

\subsection{Changes between versions 4 and 5}

Problem \ref{prob:fin} for the equations $z^2+y^2z-z+x^3+2=0$ and $z^2+y^2z+x^3+x\pm 1=0$ has been solved by ChatGPT 5.4 Pro, see \cite{grechuk2026cubic}. These equations have therefore been removed from Tables \ref{tab:H24finopen} and \ref{tab:l9openfin}. 

The equation $3x+x^2z^2+2y^2z+1=0$ has integer solutions (e.g., $x=-47103378393904407$, $y=262260511590478592$, $z=-62$). Similarly, the equation $3x^2y+y^2z^2+2z-1=0$ has integer solutions (e.g., $x=26598666324717134136290869$, $y=-141$, $z=3879814237310199004275254$). Hence, these equations have been removed from Table \ref{tab:lless11openexist}. 

\subsection{Changes between versions 5 and 6}

With the help of ChatGPT 5.4, Eugene Go found an integer solution to $1-x+x^3+x^2y^2+z+z^2=0$:
\[
(x,y,z)=(-1398651019,153,52307072551909).
\]
This had been the smallest equation in $n\geq 3$ variables for which Problem \ref{prob:main} was open; it has now been removed from Table \ref{tab:H34unsolved}.

With the help of ChatGPT 5.4, Eugene Go found an integer solution to $2+2x+x^3-y^2-xy^2+xz^2=0$:
\[
(x,y,z)=(-252123662,432516060,351431075).
\]
This had been the smallest cubic equation for which Problem \ref{prob:main} was open. ChatGPT 5.4 also found solutions to $-2x+x^3+y^2-xy^2-xz^2+3=0$ and $-2x+x^2+x^3+y+y^3+yz^2-1=0$, respectively:
\[
\begin{aligned}
(x,y,z)&=(-63017373,57059718,26746548),\\
(x,y,z)&=(-115648481,393149,1983495873).
\end{aligned}
\]
All three equations have therefore been removed.

The equation $7x^3+2y^3=3z^2+1$ was the shortest cubic equation for which Problem \ref{prob:main} was open. We have since proved that it has no integer solutions \cite{shortest}, so it has been removed. The then-new shortest cubic equations were added to Table \ref{tab:cubshortest}. The equations not appearing in \cite{mainbook} are marked in red.

\subsection{Changes between versions 6 and 7}

The solution sets of the equations $x^2+y^2+zt\pm 1=0$ are finite unions of polynomial families \cite{grechuk2026polynomial2}. For the equations $x^2y+zt+1=0$ and $x^2y+z^2\pm1=0$, the corresponding polynomial parametrizations and completeness proofs are given in \cite[Sections 8.2, 8.4 and 8.5]{grechuk2026polynomial}. These equations have therefore been removed from Table \ref{tab:H13open}.

Ashleigh Ratcliffe solved the equations $x^4-y^3+xy+x=0$, $x^4+y^3+xy-y=0$, $x^4+y^3+xy+y=0$, and $x^4+y^3+xy+x=0$ by elementary methods; see \cite[Appendix A.1]{ratcliffe2024}. These equations have been removed from Table \ref{tab:H31Open}.

Ashleigh Ratcliffe proved that both equations $x^2 y^2 + x z^2 + y \pm 1=0$ have integer solutions for every $x$ of the form $x=-(r^2+1)=-(s^2+1)/2$, where $r,s$ are integers satisfying $s^2-2r^2=1$; see \cite[Theorem A.5]{ratcliffe2024}. This solves Problem \ref{prob:large} for these equations; therefore, they have been removed from Table \ref{tab:l8openlarge}. 

\subsection{Changes between versions 7 and 8}

For the equations $y^2+z^2=x^3\pm 1$, ChatGPT 5.6 Sol Ultra produced a reasonably explicit description of all integer solutions.\footnote{See \url{https://mathoverflow.net/questions/415173/}} At that stage, only the problem of describing their rational solutions remained open.

We now regard the description in \cite[Section 4.4.1]{mainbook} of the integer solutions to $x^3+y^3+z^3=1$ as acceptable, so this equation has been removed. Likewise, $x^2y^2+xz^2+z=0$, $x^2y^2+xz^2+y=0$, and $x^3y^2=z^2\pm 1$ have been removed because the descriptions of their solution sets in \cite[Section 4.3.5]{mainbook} are now classified as ``reasonable''.

ChatGPT 5.6 Sol Ultra observed\footnote{See \url{https://mathoverflow.net/questions/484736/}} that if $(x,y,z)$ is a primitive integer solution of $x^k+2y^3+z^3=0$ for $k\in\{4,5\}$, then the integers $A=2y^3-z^3$ and $B=2yz$ satisfy $A^2+B^3=x^{2k}$. However, the resulting equations have already been solved in the literature, see \cite{bruin2003, brown2012}. Hence, the equations $x^4+2y^3+z^3=0$ and $x^5+2y^3+z^3=0$ have been removed from Section \ref{sec:allsol}.

ChatGPT 5.6 Sol Ultra observed\footnote{See \url{https://mathoverflow.net/questions/449781/}} that $y^2+x^2y+z^2x+1=0$ has solutions
\[
x=2u^2+2u+1,\qquad y=-(2u^4+2u^3+u^2-2u-1),\qquad z=(u^2+u+1)v,
\]
for every pair of integers $(u,v)$ satisfying $v^2=2u^2+2u-3$. This solves Problem \ref{prob:large} for the equation, which has therefore been removed from Tables \ref{tab:H23openlarge} and \ref{tab:l8openlarge}. Similarly, $y^2+x^2y+z^2x+2=0$ has solutions
\[
x=-2a^2,\qquad y=2(a-1)a^4+2b,\qquad z=2ba^3+1,
\]
for every pair of integers $(a,b)$ satisfying $a^2-2b^2=1$, so it has been removed from Table \ref{tab:H23openlarge}. Finally, $z^2-xy^2-x^3+2=0$ has solutions
\[
(x,y,z)=\bigl(81u^2+30u+3,\ y,\ (27u+5)(27u^2+10u+2)\bigr)
\]
for every pair of integers $(u,y)$ satisfying $3y^2=1404u^2+520u+75$; it too has been removed from Table \ref{tab:H23openlarge}.

With the help of ChatGPT 5.6 Sol Ultra, we have proved that each of the following equations has no integer solutions other than $(0,0,0,0)$:
$$
	x^3+2 x^2 y+y^3+y^2 z+x z^2+z^3+t^2 x-t^2 y-t^2 z+t y z+t^3=0,
$$
$$
	2 x^3+x^2 y+x y^2+y^3+y^2 z-2 x z^2-z^3+t^2 y+t z^2-t^3 =0,
$$
$$
	2 x^3+2 x^2 y+2 y^3+y^2 z-2 x z^2-z^3+t z^2-t^3=0.
$$
See \cite{grechuk2026long} for the proofs. 
This solves Problem \ref{prob:finhom} for these equations, so they have been removed from Section \ref{sec:fin}. The equation $8x^4+7y^4+2z^4-t^4=0$ likewise has no nonzero integer solutions \cite{grechuk2026long} and has been removed from Table \ref{tab:H288hom44}. In addition, all five equations in Table 9 of version 7 ($14x^5+5y^5+3z^5=0$, $17x^5+6y^5+3z^5=0$, $14x^5+13y^5+2z^5=0$, $21x^5+5y^5+4z^5=0$, and $17x^5+9y^5+5z^5=0$) have no nonzero integer solutions \cite{grechuk2026long}. Consequently, that table has been removed and the subsequent tables have been renumbered. Finally, $4x^5+4y^5+11z^5=0$ has no integer solutions with $z\neq 0$ \cite{grechuk2026long}, so it has been removed from Section \ref{sec:large}.

Ashleigh Ratcliffe proved that the equation $x^3y^2+z^2+y+1=0$ admits integer solutions in which all coordinates are arbitrarily large in absolute value; see \cite[Theorem A.7]{ratcliffe2024}. This solves Problem \ref{prob:large} for this equation, which has therefore been removed from Table \ref{tab:l8openlarge}.

ChatGPT 5.6 Sol Ultra observed\footnote{See \url{https://mathoverflow.net/questions/450943/}} that $y(x^3-z^2)=z+1$ has solutions
\[
x=9u^4+8u,\qquad y=-\frac{9(3u^3+1)}{64},\qquad z=27u^6+36u^3+8,
\]
which are integral whenever $u\equiv 13\pmod {64}$. In addition, for every pair of integers $(a,b)$ satisfying $a^2-5b^2=-4$, the formulas
\[
x=\frac{a^2+12a+16}{5},\qquad y=-\frac{25(a+1)}{1728},\qquad z=\frac{b(a^2+18a+76)}{5}
\]
give a rational solution of $y(x^3-z^2)=x-1$; suitable congruence conditions on $a$ and $b$ make these values integral. This solves Problem \ref{prob:large} for both equations, which have therefore been removed from Table \ref{tab:l8openlarge}.

The equation $y^2+x^3y+z^2+1=0$ has infinitely many integer solutions \cite[Corollary 3.1]{forms}.\footnote{See also \url{https://mathoverflow.net/questions/454040/}.} This also solves Problem \ref{prob:large} for this equation: fixing $x$ bounds $y$ and $z$ by completing the square; fixing $z$ leaves finitely many possibilities because $y\mid z^2+1$; and fixing a nonzero $y$ gives a nonsingular genus-$1$ curve, with finitely many integer points by Siegel's theorem. The case $y=0$ is impossible. Thus, only finitely many solutions have any one coordinate bounded in absolute value. The equation has therefore been removed from Tables \ref{tab:l8openlarge} and \ref{tab:l9openfin}. The same method proves infinitude for $y^2+x^3y+z^2-2=0$, $y^2+x^3y+z^2+z-1=0$, $y^2+x^3y+z^2+z+1=0$, and $y^2+x^3y+y+z^2+1=0$ \cite[Corollary 3.2]{forms}; these have therefore been removed from Table \ref{tab:l9openfin}.

\subsection{Changes between versions 8 and 9}

Jamal Agbanwa, with the help of ChatGPT 5.6 Sol Ultra, found infinitely many integer solutions to $z^2+y^2z+2x^3+1=0$. For every odd integer $r$, set
\[
\begin{aligned}
x&=-\frac{r^8-8r^6+22r^4-22r^2+5}{2},\\
y&=r^3-3r,\\
z&=\frac{r^{12}-12r^{10}+57r^8-134r^6+159r^4-84r^2+11}{2}.
\end{aligned}
\]
These formulas give an integer solution for every such $r$. Hence, the equation has been removed from Table \ref{tab:l9openfin}.

Jamal Agbanwa, with the help of ChatGPT 5.6 Sol Ultra, also found infinitely many integer solutions to $x^3+x^2y^2+z^2+1=0$. Let
\[
Q(t)=465577896+1881226506t+1900333542t^2.
\]
The Pell-type equation $Y^2=Q(t)$ has infinitely many integer solutions $(t,Y)$; for example, $Q(-659)=28706088^2$, and the standard Pell recurrence generates infinitely many further solutions. For each such pair, set
\[
\begin{aligned}
x&=-577630421-2333988934t-2357694538t^2,\\
y&=Y,\\
z&=50421655389668t^3+74872031013786t^2+37059612550518t+6114499038718.
\end{aligned}
\]
Direct substitution gives $x^3+x^2y^2+z^2+1=0$. Hence, the equation has been removed from Table \ref{tab:l9openfin}.

Using ChatGPT 5.6 Sol Ultra, Eugene Go showed that $z^2-xy^2-x^3-2=0$ has integer solutions in which all variables are arbitrarily large in absolute value. Set $d=8353$. The generalized Pell equation
\[
Y^2-2(d^3+1)r^2=d^6-4
\]
has an integer solution with $d^2\mid Y$. Standard Pell theory, applied with a suitable power of the fundamental unit modulo $d^2$, then gives infinitely many solutions with this divisibility property. Writing $Y=d^2y$ gives infinitely many integer pairs $(y,r)$ satisfying
\[
d^4y^2-2(d^3+1)r^2=d^6-4.
\]
For each such pair, define
\[
w=\frac{dy^2-2r^2-d^3}{2},\qquad x=dw,\qquad z=r(w+1).
\]
Here $w$ is an integer. The Pell equation implies $r^2=d^3w+2$, and direct calculation gives
\[
z^2-xy^2-x^3-2=(w^2+1)(r^2-d^3w-2)=0.
\]
Hence, the equation has been removed from Table \ref{tab:H23openlarge}.

Using ChatGPT 5.6 Sol Ultra, Eugene Go also showed that $y^2+x^2y+z^2x-2=0$ has integer solutions in which all variables are arbitrarily large in absolute value. For integers $p,q$ satisfying $p^2-6pq+q^2=1$, define
\[
\begin{aligned}
T&=128p^8q^2-4p^2+1,\\
x&=2-64p^6q^2,\\
y&=8p^3qT,\\
z&=(q-p)T-2p.
\end{aligned}
\]
These formulas give an integer solution, and the Pell-type equation for $(p,q)$ supplies solutions of arbitrarily large size. Hence, the original equation has been removed from Table \ref{tab:H23openlarge}; in its place, equations of size $H=23$ for which Problem \ref{prob:large} remains open have been added to the table.

Using ChatGPT 5.6 Sol Ultra, Eugene Go proved that the equations $x^4+x^2y^2+y^3z-yz^3=0$, $x^4-y^4+x^2yz+yz^3=0$, $x^4+y^4+x^2yz-yz^3=0$, $x^4-x^2y^2+y^3z+yz^3=0$, and $x^4+x^2y^2+y^3z+z^4=0$ have no integer solutions with $xyz\neq 0$. See \cite{grechuk2026long} for the proofs. These equations have therefore been removed from Table \ref{tab:H64genus3}.

Table \ref{tab:l9openlargecub} has been added to list the shortest cubic equations for which Problem \ref{prob:large} remains open. The equation $2x^3+2y^3+2z^3=xyz+1$ has been added as another shortest symmetric (and cyclic) equation for which Problem \ref{prob:main} remains open; the next shortest is $x^3+y^3+z^3=114$. The equation $3x^4+y^4=z^4+t^4$ has been added as both the smallest and the shortest equation of the form \eqref{eq:homdiand} for which Problem \ref{prob:homlarge} remains open. Likewise, $x^4+y^4+z^4=51t^4$ has been added as the shortest equation of the form \eqref{eq:homdiand} for which Problem \ref{prob:finhom} remains open.

A paragraph describing the additional cross-check of the stated minimality claims has been added to the introduction.

\subsection{Changes between versions 9 and 10}

With the help of GPT-6 Astra Pro, we proved that the solution sets of the equations $xyz+t^2-1=0$ and $xyz+ts+1=0$ are finite unions of polynomial families \cite{grechuk2026polynomial}. These equations have therefore been removed from Table \ref{tab:H13open}.

With the help of GPT-6 Astra Pro, we determined all integer solutions to $x^4+y^3+x-2y=0$ and $x^4+y^3+x-2y+1=0$, see \cite[Proposition 12]{grechuk2026}. Using GPT-6 Astra Pro and Magma, Eugene Go proved that the only integer solution to $x^4+y^3+xy+x+1=0$ is $(x,y)=(0,-1)$; see \cite[Theorem 3.14.1]{grechuk2026long}. All three equations have therefore been removed from Table \ref{tab:H31Open}.

Using GPT-6 Astra Pro, Avi Eisenberg determined all primitive integer solutions to the seven equations in Table 4 of version 9: $x^4+3y^3+z^3=0$, $x^4+3y^3+2z^3=0$, $2x^4-y^4+z^3=0$, $x^5+y^4-2z^2=0$, $x^4+4y^3+z^3=0$, $x^5+y^4-3z^2=0$, and $x^5+y^4+3z^2=0$, see  \cite{grechuk2026long}. Here primitive means $\gcd(x,y,z)=1$. That table and the accompanying generalized Fermat discussion have therefore been removed from Section \ref{sec:allsol}.

Using Claude Opus 5.5, Avi Eisenberg proved that $x^4+x^3y+xy^3-z^4=0$ has no integer solutions with $xyz\neq 0$, see \cite{grechuk2026long}. This equation has therefore been removed from Table \ref{tab:H64genus3}.

Using GPT-6 Astra Ultra, Eugene Go found the integer solution
\[
(x,y,z,t)=(7337397,6619385,8035045,12352344)
\]
to $7x^4+5y^4+4z^4-2t^4=0$. This equation has therefore been removed from Table \ref{tab:H288hom44}.

Using GPT-6 Astra Pro, Eugene Go proved that all 30 equations in Table 9 of version 9 have no integer solutions. These are
\begingroup
\[
\begin{array}{ll}
y^3+xy=x^4+4, & y^3+xy=x^4+x+2,\\
y^3+y=x^4+x+4, & y^3-y=x^4+2x-2,\\
y^3+y^2+xy=x^4+1, & y^3+xy=x^4+x+3,\\
y^3+xy=x^4-x-3, & y^3+xy=x^4+x-3,\\
y^3+xy+y=x^4-x+1, & y^3+y^2=x^4+x+4,\\
y^3-y^2=x^4+2x+2, & y^3+y^2=x^4+2x-2,\\
y^3-y^2+xy=x^4+2, & y^3+y=x^4+x+6,\\
y^3+y=x^4+x-6, & y^3+2y=x^4+x+4,\\
y^3+y=x^4+2x+4, & y^3+y=x^4+2x-4,\\
y^3+xy=x^4+x+4, & y^3+xy=x^4-x-4,\\
y^3+xy=x^4+x^2+2, & y^3+y^2+y=x^4+x+2,\\
y^3+y^2-y=x^4+x+2, & y^3-y=x^4-x^2+x+2,\\
y^3-y=x^4-x^2+x-2, & y^3+xy+y=x^4+4,\\
y^3+xy-y=x^4+4, & y^3-y^2+y=x^4+x+2,\\
y^3+xy-y=x^4+x+2, & y^3+xy-y=x^4-4.
\end{array}
\]
\endgroup
The proofs can be found in \cite{grechuk2026}. That table has therefore been removed. Table \ref{tab:H34unsolved} now lists the smallest open equations without restricting the number of variables. Eugene Go also proved that $2y^3+xy+x^4+1=0$ has no integer solutions, see \cite{grechuk2026}, so this equation has been removed from Table \ref{tab:lless11openexist}.

Using GPT-6 Astra Ultra, Eugene Go (and, independently, Avi Eisenberg using Claude Opus 5.5) observed that the substitutions
\begingroup
\[
x=-\frac{2r^2+1}{d^2},\qquad
z=r\left(\frac{2r^2+1}{d^3}+1\right)
\]
\endgroup
transform the left-hand side of $2z^2+xy^2+x^3+1=0$ into
\begingroup
\[
\frac{2r^2+1}{d^6}
\left(d^6-1+(4d^3-2)r^2-d^4y^2\right).
\]
\endgroup
It therefore suffices to find an integer $d>1$ for which
\begingroup
\[
d^4y^2-(4d^3-2)r^2=d^6-1
\]
\endgroup
has an integer solution $(y,r)$. Direct search shows that $d=2689$ works. To obtain infinitely many such pairs, write $Y=d^2y$. Since $4d^3-2>0$ and $4d^3-2\equiv 2\pmod 4$, this coefficient is not a square, and $Y^2-(4d^3-2)r^2=d^6-1$ is a generalized Pell equation with nonzero right-hand side. A suitable power of a Pell unit preserves $d^2\mid Y$ and generates infinitely many solutions with $|r|,|y|\to\infty$. Reducing the displayed equation modulo $d^3$ gives $2r^2+1\equiv 0\pmod{d^3}$, so the formulas for $x$ and $z$ give integers. Along this sequence, $|x|$ and $|z|$ also tend to infinity. This solves Problem \ref{prob:large}, so $2z^2+xy^2+x^3+1=0$ has been removed from Table \ref{tab:l9openlargecub}.

Using Claude Opus 5.5, Avi Eisenberg found the integer solution
\[
(x,y,z)=(16852584801875,20167801990609,-28496089784237)
\]
to $y^2z+yz^2=x^3+x^2+3x-1$. This was previously the smallest open cubic equation for Problem \ref{prob:main}. Equation $x^2 y+x y^2+x z^2+2x^2-z^2+3=0$ was added alongside \eqref{eq:H39sym}; these were then the two smallest open cubic equations, both of size $H=39$.

Using Claude Opus 5.5, Avi Eisenberg also found that the equation $2x^3+3xy^2+z^3+z^2+1=0$ has the integer solution
\[
(x,y,z)=(-15985616853,84362449981,70438064692),
\]
and $y^2z+yz^2=3x^3+x^2+x-1$ has the integer solution
\[
(x,y,z)=(-783954692511,1797526169071,-838276125548).
\]
These two equations, together with $y^2z+yz^2=x^3+x^2+3x-1$, have therefore been removed from Table \ref{tab:cubshortest}. Further, the equation $x^3+2x+y^3+2y+z^3+2z=xyz+1$ has the integer solution
\[
(x,y,z)=(723809206820,-631758385864,-245500254229),
\]
so it has been removed from Section \ref{sec:existence}, while the equation $2x^3+2y^3+2z^3=xyz+1$ has the integer solution
\[
(x,y,z)=(784464819555,-706181253695,-335317485833),
\]
so it has been removed from Section \ref{sec:shortopen}. The corresponding smallest and shortest open symmetric and cyclic equation statements have been updated.

Following the deletion of Tables 4 and 9 of version 9, subsequent tables have been renumbered automatically. References to the deleted tables in earlier change summaries have been replaced by section references. 

\subsection{Changes between versions 10 and 11}

Equations \eqref{eq:zy2pm1x2p2} are added as examples of equations with uncertain status. 

Using Claude Opus 5.5, Eugene Go described all rational solutions to $y^2+z^2=x^3+1$, $y^2+z^2=x^3-1$, and $y^2-x^2y+z^2+1=0$.\footnote{See Eugene Go's answer of 26 October 2026 at \url{https://mathoverflow.net/a/516407}. The plus-sign form of the third equation in the original question is equivalent under $y\mapsto-y$.} These rational points correspond to primitive integer solutions, up to simultaneous sign and with $t\neq0$, of the respective homogenizations
\[
y^2t+z^2t-x^3-t^3=0,\qquad
y^2t+z^2t-x^3+t^3=0,
\]
\[
y^2t-x^2y+z^2t+t^3=0.
\]
All integer solutions with $t\neq0$ are obtained by integer scaling of primitive representatives. For $t=0$, the first two equations give $(x,y,z,t)=(0,y,z,0)$, while the third gives $x^2y=0$, with $z$ arbitrary. Thus the affine rational descriptions, together with these explicit solutions at infinity, describe all integer solutions of the homogenizations. These equations have therefore been removed from Section \ref{sec:allsol}.

Using GPT-6 Astra Pro, Eugene Go determined all integer solutions to $x^4-x^2+xy+y^3=0$, which was listed in Section \ref{sec:allsol} as the smallest open two-variable equation whose normalization has genus $2$; see \cite[Theorem 3.15.1]{grechuk2026long}. This equation has therefore been removed.

Using GPT-6 Astra Pro, all integer solutions to the following equations have been determined:
\[
\begin{array}{ll}
	x^4-y^3-x+y=0, & x^4-y^3+x-y=0,\\
	x^4+y^3+xy+1=0, & x^4+y^3+2x-y=0,\\
	x^4+y^3+2x+y=0, & x^4+y^3+y^2+x=0,\\
	x^4+y^3-y^2+x=0, & x^4+y^3+2x-y+1=0,\\
	x^4+y^3+2x+y+1=0, & x^4+y^3+2x+y-1=0,\\
	x^4+y^3+xy-y-1=0, & x^4-y^3+xy+x-1=0,\\
	x^4+y^3+xy+y-1=0.
\end{array}
\]
See \cite{grechuk2026long} for the proofs. These equations have therefore been removed from Table \ref{tab:H31Open}. The three corresponding length-$9$ entries have also been removed from Table \ref{tab:l9openfin}.

Using Claude Opus 5.5, Avi Eisenberg proved, initially assuming the generalized Riemann hypothesis (GRH), that $x^3z+y^3x+y^2z^2-z^3y=0$ and $x^4-y^3z+xyz^2+xz^3=0$ have no integer solutions satisfying $xyz\neq0$. Eugene Go subsequently removed the GRH assumption. The unconditional results, including complete descriptions of the solutions with $xyz=0$, are proved in \cite{grechuk2026long}. These equations have therefore been removed from Table \ref{tab:H64genus3}.

Using GPT-6 Astra Ultra, Eugene Go found the integer solution
$$
	(x, y, z, t) = (32519425, 30486957, 21147587, 63400357)
$$
to $11x^4+4y^4+2z^4-t^4=0$, as well as the integer solution
$$
	(x, y, z, t) = (60428235, 11074181, 51780065, 86883262)
$$
to the equation $8x^4 + 7y^4 + z^4 - 2t^4 = 0$. Hence, these equations have been removed from Table \ref{tab:H288hom44}.

Using Claude Opus 5.5, Jamal Agbanwa found an integer solution
$$
	(x,y,z) = (-2450412151997178, 92792217037094, 467727250324161)
$$
to the equation $x^2 y+x y^2+x z^2+2x^2-z^2+3=0$, which was, together with \eqref{eq:H39sym}, one of the smallest open cubic equations for Problem \ref{prob:main}. This equation has now been removed from Section \ref{sec:existence}. 

Table \ref{tab:version_progress} has been added to summarize progress on the equations from the original book. Here ``solved'' includes equations whose existing solution descriptions have been reclassified as acceptable.

\begin{table}
	\begin{center}
		\begin{tabular}{ |c|c|c|c|c| }
			\hline
			Version & Date & Newly solved & Total solved & Remaining open \\
			\hline\hline
			v1  & 12 Apr 2024 & --- & 0 & 205 \\
			\hline
			v2  & 15 Oct 2024 & 13 & 13 & 192 \\
			\hline
			v3  & 26 Dec 2024 & 4 & 17 & 188 \\
			\hline
			v4  & 23 Jan 2026 & 4 & 21 & 184 \\
			\hline
			v5  & 8 Mar 2026 & 5 & 26 & 179 \\
			\hline
			v6  & 1 Apr 2026 & 5 & 31 & 174 \\
			\hline
			v7  & 22 May 2026 & 11 & 42 & 163 \\
			\hline
			v8  & 13 Jul 2026 & 27 & 69 & 136 \\
			\hline
			v9  & 30 Aug 2026 & 9 & 78 & 127 \\
			\hline
			v10 & 1 Oct 2026 & 46 & 124 & 81 \\
			\hline
			v11 & This version & 21 & 145 & 60 \\
			\hline
		\end{tabular}
		\caption{\label{tab:version_progress} Numbers of equations from the original book newly solved and remaining open in successive versions of this document. Equations reducible to each other or appearing under multiple problems are counted once; equations marked in red are excluded.}
	\end{center}
\end{table}

\bibliography{open_main_refs}
\bibliographystyle{plain}

\end{document}